\documentclass[12pt,cs4size,a4paper,fancyhdr,fntef]{article}
\ifx\pdfoutput\undefined 
  \usepackage[dvipdfm,CJKbookmarks]{hyperref}
\else
  \usepackage[]{hyperref}
\fi
\usepackage{shortvrb,ulem,makeidx} 
\usepackage{fancyhdr}
\usepackage{multicol}
\usepackage{array,color}
\usepackage{amssymb,amsthm}
\usepackage{amsfonts}
\usepackage[reqno]{amsmath}
\usepackage{extarrows}
\usepackage{graphicx}
\usepackage{makecell,rotating}
\usepackage{multirow}
\usepackage{url}
\usepackage{indentfirst}
\usepackage{subfigure}

 \MakeShortVerb{\|}

 \numberwithin{equation}{section}

\newtheorem{defi}{DEFINITION}[section]
\newtheorem{theo}[defi]{THEOREM}
\newtheorem{prop}[defi]{PROPOSITION}
\newtheorem{lemm}[defi]{LEMMA}
\newtheorem{coro}[defi]{COROLLARY}

\newtheorem{rem}[defi]{REMARK}

\makeatletter
\def\@biblabel#1{#1}
\makeatother

\begin{document}

\centerline{\large\bf Some apriori estimates of G-BSDEs and }
\centerline{\large\bf the G-martingale representation for a special case}

\vskip 0.6cm

\centerline{\bf FAN Yulian }

\vskip0.2cm

\centerline{School of Science, North China University of Technology, Beijing 100144, China}
\centerline{Email: fanyl@ncut.edu.cn}

\date{}

\vskip 0.3cm

\noindent{\bf Abstract} This paper presents the integral(or differential) form of G-BSDEs, gives some kind of apriori estimates of their solutions, and under a very strong condition, proves the G-martingale representation theorem, and the existence and uniqueness theorem of G-BSDEs.

\vskip0.5cm


\baselineskip13.9pt \abovedisplayskip=4pt plus 1pt minus 1pt
\abovedisplayshortskip=4pt plus 1pt minus 1pt \belowdisplayskip=4pt
plus 1pt minus 1pt \belowdisplayshortskip=4pt plus 1pt minus 1pt \
\vskip -1cm



\vskip0.6cm

\section{Introduction}

Backward stochastic differential equations(BSDEs) was first introduced by Pardoux and Peng [8] in the classical
linear expectation case. Since then on, a lot of works have been devoted to study BSDE theory and its applications.

Based on BSDE theory, Peng [9]  introduced the notion of g-expectation and
conditional g-expectation which is the first dynamically consistent nonlinear expectation.

Then Peng [10] introduced the notion of G-expectation, which is a more general dynamically consistent nonlinear expectation,  and the concept of G-Brownian motion, and then established the related stochastic calculus. The theory of G-expectation is intrinsic in the
sense that it is not based on a given (linear) probability space, and it takes the probability uncertainty into consideration. Drift uncertainty and volatility uncertainty are two typical situations of probability uncertainty.
G-Brownian motion has a very rich and interesting
new structure which non-trivially generalizes the classical Brownian motion. G-expectation theory has developed rapidly since the initial paper Peng [10].  
Peng [11, 13, 14]  
studied the central limit theorem under
sublinear expectations and obtained that the limit distribution exists and is just the
G-normal distribution.
Peng [12]  
and
Peng [15] 
systematically developed the stochastic calculus under G-expectation.
Xu and Zhang [19]  
studied
the It\^o's integral with respect to G-martingales and the L\'{e}vy characterization of G-Brownian
motion. Gao [3] 
studied the path properties of the solutions of G-SDEs.
Hu and Peng [5] 
studied G-L\'{e}vy processes. Li and Peng [7] 
studied the It\^o's integral without the condition of quasi-continuous and on stopping time interval, and generalized the It\^o's formula to general $C^{1,2}$-function. Denis, Hu and Peng [2] 
and Hu and Peng [5] 
studied the representation theorem of G-expectation
and its application to G-Brownian motion paths. Soner, Touzi and Zhang [17], 
Song [18], 
Hu and Peng [6], 
and Peng, Song and Zhang [16] 
studied the G-martingale representation theorem.

G-Brownian motion has independent increments with identical
G-normal distributions which means it can characterize the volatility uncertainty. And a very interesting new phenomenon  is that its quadratic process generally is not a deterministic process but a stochastic process which also has independent increments
with identically maximal distributions. So the stochastic differential equations driven by G-Brownian motion(G-SDEs)
of the following form
\begin{equation*}\label{SDE}
X_{s}=b(s,X_{s})ds+h(s,X_{s}):d\langle B\rangle_{s}+\theta^{*}(s,X_{s})dB_{s}, t\in [0,T]
\end{equation*}
will carry the characteristic of both mean uncertainty and volatility uncertainty.
Peng [10] 
proved the existence and uniqueness of the solutions of such G-SDEs. Gao [3] 
gave some moment estimates and
H\"{o}lder continuity results of the solution of G-SDEs. However the corresponding problems for backward
stochastic differential equations are not completely solved. Peng [15] 
give partial results to this direction, i.e.,
the following type of G-BSDE:
\begin{equation*}\label{G-BSDE expectation form}
Y_{t}=\mathbb{E}[\xi+\int_{t}^{T}f(s,Y_{s})ds+\int_{t}^{T}h(s,Y_{s}):d\langle B\rangle_{s}\mid\Omega_{t}], t\in [0,T].
\end{equation*}
has a unique solution if the coefficients $f,h$ satisfies Lipschitz condition.
Hu, et al (2012) prove the existence and uniqueness result of the following G-BSDE
\begin{equation*}\label{G-BSDE of Hu et al 2012}
Y_{t}=\xi+\int_{t}^{T}f(s,Y_{s},Z_{s})ds+\int_{t}^{T}g(s,Y_{s},Z_{s}):d\langle B\rangle_{s}-\int_{t}^{T}Z^{*}_{s} dB_{s}-(K_{T}-K_{t})
\end{equation*}
 by applying the partition of unity theorem to construct
a new type of Galerkin approximation.

The aim of this paper is to give some kind of apriori estimates of G-BSDE
\begin{equation*}\label{G-BSDE integral form}
Y_{t}=\xi+\int_{t}^{T}f(s,Y_{s},Z_{s},\eta_{s})ds+\int_{t}^{T}g(s,Y_{s},Z_{s},\eta_{s}):d\langle B\rangle_{s}-\int_{t}^{T}Z^{*}_{s} dB_{s}+\int_{t}^{T}G(\eta_{s})ds-\frac{1}{2}\int_{t}^{T}\eta_{s}:d\langle B\rangle_{s}.
\end{equation*}
or, equivalently, the differential form
\begin{equation*}\label{G-BSDE differential form}
-dY_{t}=f(t,Y_{t},Z_{t},\eta_{t})dt+g(t,Y_{t},Z_{t},\eta_{t}):d\langle B\rangle_{t}-Z^{*}_{t} dB_{t}+G(\eta_{t})dt-\frac{1}{2}\eta_{t}:d\langle B\rangle_{t}, Y_{T}=\xi,
\end{equation*}
and under a very strong condition, we get the G-martingale representation theorem, and the existence and uniqueness of the solution $(Y,Z,\eta)$.

The rest of this paper is organized as follows. In section 2, we introduce the notations and definitions. In section 3, we give some kind of apriori estimates. In section 4, under a very strong condition, we get the G-martingale representation theorem and the existence and uniqueness results of the solutions of G-BSDEs.

\section{Preliminaries}

For any $n\times d$ dimensional matrices $\gamma, \tilde{\gamma}$, define
\begin{equation*}
\gamma:\tilde{\gamma}:=tr(\gamma^{*}\tilde{\gamma}), \vert\gamma\vert:=\sqrt{\gamma:\gamma},
\end{equation*}
where $\gamma^{*}$ denotes the transpose of $\gamma$.

For a dimension $d$, let
$\mathbb{R}^{d}$, $\mathbb{S}^{d}$, and $\mathbb{D}^{d}$ denote the sets of $d$-dimensional  column vectors, $d\times d$-symmetric matrices, and $d\times d$-diagonal matrices respectively. For $\sigma_{1},\sigma_{2}\in\mathbb{S}^{d}$, $\sigma_{1}\leq\sigma_{2}$(resp. $\sigma_{1}<\sigma_{2}$) means that $\sigma_{2}-\sigma_{1}$ is nonnegative(resp. positive) definite, and we
denote by $[\sigma_{1},\sigma_{2}]$ the set of $\sigma\in\mathbb{S}^{d}$ satisfying $\sigma_{1}\leq\sigma\leq\sigma_{2}$. For $\xi^{1},\xi^{2}\in \mathbb{R}^{d}$, $\xi^{1}\leq\xi^{2}$ means that each
element of $\xi^{1}$ is less or equal to that of $\xi^{2}$, that is $\xi^{1}_{i}\leq \xi^{2}_{i}, i=1,\ldots,d$. We use $\mathbf{0}$ to denote the $d$-dimensional zero vector or zero matrix, and $I_{d}$ the $d\times d$ identity matrix.
And for $\gamma, \tilde{\gamma}\in \mathbb{S}^{d}$, we have
\begin{equation*}
\vert\gamma:\tilde{\gamma}\vert\leq \vert \gamma\vert\vert \tilde{\gamma}\vert, \text{and} -\gamma\leq\tilde{\gamma}\leq\gamma \text{ implies that } \vert\tilde{\gamma}\vert\leq \vert\gamma\vert.
\end{equation*}

Let $\mathbb{R}^{n\times d\times d}$ denote all $\eta=(\eta^{1},\ldots,\eta^{n})^{*}$ with $\eta^{i}, i=1,\ldots,n$ being $d\times d$ matrices. When $\eta^{i}, i=1,\ldots,n$ are symmetric matrices, we use $\mathbb{S}^{n\times d\times d}$ instead of  $\mathbb{R}^{n\times d\times d}$, when  $\eta^{i}=diag(\eta^{i;1},\ldots,\eta^{i;d}), i=1,\ldots,n$  are diagonal matrices, we use $\mathbb{D}^{n\times d \times d}$ instead of $\mathbb{S}^{n\times d\times d}$, and when $n=1$, we use $\mathbb{D}^{d}$ instead of $\mathbb{D}^{1\times d \times d}$. For any symmetric matrix $\gamma$, define
\begin{equation*}
\eta:\gamma=(\eta^{1}:\gamma, \ldots, \eta^{n}:\gamma)^{*},
\end{equation*}
with $\eta^{i}:\gamma=tr((\eta^{i})^{*}\gamma), i=1,\ldots,n$.
Now we define an operator $\cdot$, such that
\begin{equation*}
\begin{split}
\eta\cdot\theta&=\sum_{i=1}^{n}(\eta^{i})^{*}\theta^{i}, \text{ if } \eta,\theta\in \mathbb{R}^{n\times d \times d},\\
\xi\cdot\eta&=\eta\cdot\xi=\sum_{i=1}^{n}\xi^{i}\eta^{i}, \text{ if } \xi\in \mathbb{R}^{n}, \eta\in \mathbb{R}^{n\times d \times d},\\
\end{split}
\end{equation*}
and
\begin{equation*}
\xi\cdot\eta:\gamma=\sum_{i=1}^{n}\xi^{i}\eta^{i}:\gamma, \text{ if } \xi\in \mathbb{R}^{n}, \eta\in \mathbb{R}^{n\times d \times d}, \gamma\in \mathbb{S}^{d}.
\end{equation*}
Define
\begin{equation*}
\vert\eta\vert=\sqrt{\sum_{i=1}^{n}\eta^{i}:\eta^{i}},\text{ if } \eta\in \mathbb{R}^{n\times d \times d},
\end{equation*}
and $G(\eta)=(G(\eta^{1}),\ldots,G(\eta^{n}))^{*}$ where
\begin{equation*}
G(\eta^{i}):=\frac{1}{2}\sup_{\sigma^{2}\in[\underline{\sigma}^{2},\bar{\sigma}^{2}]}(\sigma^{2}:\eta^{i}).
\end{equation*}

Let $\Omega=C([0,\infty],\mathbb{R}^{d})$, $\mathcal{F}=\mathcal{B}(\Omega)$. G-expectation $\mathbb{E}$ is
a sublinear expectation on the canonical space $\Omega$ such that the canonical process $B$ is a G-Brownian motion. We assume the increment
$B_{t+s}-B_{t}$ is $N(\{0\}\times\Sigma s)$-distributed, for each $t,s\geq0$, where $\Sigma$ is a bounded, convex and closed subset of $\mathbb{S}^{d}$.

Define
\begin{equation*}
G(A):=\frac{1}{2}\mathbb{E}[\langle AB_{1}, B_{1}\rangle]=\frac{1}{2}\sup_{\sigma^{2}\in\Sigma} tr[A\sigma^{2}], A\in \mathbb{S}^{d}.
\end{equation*}

We assume $B_{t}=(B^{1}_{t},\ldots,B^{d}_{t})^{*}$ satisfies that for each fixed $t$, $B^{i}_{t},i=2,\ldots,n$ is independent from
$B^{i-1}_{t},\ldots, B^{1}_{t}$.
Then it is easy to prove that the matrices in $\Sigma$ will be diagonal matrices, i.e., any $\sigma^{2}\in\Sigma$, $\sigma^{2}=diag(\sigma^{2}_{11},\ldots,\sigma^{2}_{nn})$, with
$\sigma^{2}_{ii}\in[\underline{\sigma}^{2}_{ii}, \bar{\sigma}^{2}_{ii}]$, where
$\underline{\sigma}^{2}_{ii}=-\mathbb{E}[-B^{i}_{1}B^{i}_{1}], \bar{\sigma}^{2}_{ii}=\mathbb{E}[B^{i}_{1}B^{i}_{1}]$.
In the following, we denote  $\underline{\sigma}^{2}=diag(\underline{\sigma}^{2}_{11},\ldots,\underline{\sigma}^{2}_{nn})$, and $\bar{\sigma}^{2}=diag(\bar{\sigma}^{2}_{11},\ldots,\bar{\sigma}^{2}_{nn})$.

Assume $\tilde{B}$ is another G-Brownian motion without the independence assumption, then $\tilde{B}_{t+s}-\tilde{B}_{t}\sim N(\{0\}\times\tilde{\Sigma}s)$ where $\tilde{\Sigma}$ is a bounded, convex and closed subset of $\mathbb{S}^{d}$. Define
\begin{equation*}
\tilde{G}(A):=\frac{1}{2}\tilde{\mathbb{E}}[\langle A\tilde{B}_{1}, \tilde{B}_{1}\rangle]=\frac{1}{2}\sup_{Q\in\tilde{\Sigma}} tr[AQ], A\in \mathbb{S}^{d}.
\end{equation*}

Denote $Q=(q_{ij}), \underline{Q}=(\underline{q}_{ij}), \bar{Q}=(\bar{q}_{ij})$,
where $\underline{q}_{ij}=\inf q_{ij}=-\mathbb{E}[-\tilde{B}^{i}_{1}\tilde{B}^{j}_{1}]$,
$\bar{q}_{ij}=\sup q_{ij}=\mathbb{E}[\tilde{B}^{i}_{1}\tilde{B}^{j}_{1}]$.

Denote $X_{i}=B^{i}_{1}, i=1,\ldots,d$ and $X=(X_{1},\ldots,X_{d})^{*}$. Let $P$ be a matrix, and let
$Y=PX$ , then
\begin{equation*}
\begin{split}
&\frac{1}{2}\mathbb{E}[\langle AY, Y\rangle]=\frac{1}{2}\mathbb{E}[\langle APX,PX\rangle]=\frac{1}{2}\mathbb{E}[\langle P^{*}APX,X\rangle]\\
=&\frac{1}{2}\sup_{\sigma^{2}\in\Sigma}  tr(P^{*}AP\sigma^{2})=G(P^{*}AP).
\end{split}
\end{equation*}

We let the elements in matrix $A=(a_{ij})$ take the following values.
For given $i,j=1,\ldots,d$, we let $a_{ij}=a_{ji}=1$, $a_{kl}=0$, if $k,l=1,\ldots,n, (k,l)\neq (i,j)$ and $(j,i)$, and  then
\begin{equation}\label{ij}
\sup_{\sigma^{2}\in\Sigma}\sum_{l=1}^{n}(p_{il}p_{jl}+p_{jl}p_{il})\sigma_{ll}^{2}=2\mathbb{E}[Y_{i}Y_{j}],  i,j=1,\ldots,d.
\end{equation}
If we let the 1 in above procedure be replaced by $-1$, then we get
\begin{equation}\label{-ij}
\inf_{\sigma^{2}\in\Sigma^{2}}\sum_{l=1}^{n}(p_{il}p_{jl}+p_{jl}p_{il})\sigma_{ll}^{2}=-2\mathbb{E}[-Y_{i}Y_{j}], i,j=1,\ldots,d.
\end{equation}
So as long as there exist matrix $P$ and diagonal matrices set $\Sigma$ such that (\ref{ij}) and (\ref{-ij}) hold, we can construct random vector $\tilde{B}_{1}=Y=PX$ whose covariance matrices set is $\tilde{\Sigma}$ such that any $Q=(q_{ij})\in\tilde{\Sigma}$, $\underline{q}_{ij}\leq q_{ij}\leq\bar{q}_{ij}$ where $\underline{q}_{ij}=\inf q_{ij}=-\mathbb{E}[-\tilde{B}^{i}_{1}\tilde{B}^{j}_{1}]$,
$\bar{q}_{ij}=\sup q_{ij}=\mathbb{E}[\tilde{B}^{i}_{1}\tilde{B}^{j}_{1}]$. We denote $\underline{Q}=(\underline{q}_{ij}), \bar{Q}=(\bar{q}_{ij})$.

So in this paper we assume for each fixed $t$, $B^{i}_{t},i=2,\ldots,d$ is independent from
$B^{i-1}_{t},\ldots, B^{1}_{t}$. Hence $\Sigma$ is bounded, convex and closed subset of diagonal matrices $\mathbb{D}^{d}$.

Let $\langle B\rangle$ denote the quadratic variation of $B$  such that $BB^{*}-\langle B\rangle$ is a G-martingale. Since for each fixed $t$, $B^{i}_{t},i=2,\ldots,d$ is independent from $B^{i-1}_{t},\ldots, B^{1}_{t}$,  we can let $\langle B\rangle_{t}$  be a diagonal matrix for every $t$.

For each fixed $T\geq 0$, let
$$
L_{ip}(\Omega_{T}):=\{\varphi(B_{t_{1}},B_{t_{1}},\ldots,B_{t_{n}}):\forall
n\geq 1, t_{1},\ldots,t_{n}\in[0,T], \forall \varphi\in
C_{l,Lip}(\mathbb{R}^{d\times n}) \}.
$$
For fixed $p\geq 1$, define a norm on $L_{ip}(\Omega_{T})$
\begin{equation*}
\parallel \xi\parallel_{L^{p}_{G}(\Omega_{T})}^{p}=(\mathbb{E}[\vert \xi\vert^{p}])^{\frac{1}{p}},
\end{equation*}
and denote $L^{p}_{G}(\Omega_{T})$ be the closure of $L_{ip}(\Omega_{T})$ under the norm $\parallel \cdot\parallel_{L^{p}_{G}(\Omega_{T})}$.

$\mathbb{L}^{p}$ denote the Banach space under the norm $\parallel X\parallel_{p}:=(\mathbb{E}[\vert \xi\vert^{p}])^{\frac{1}{p}}$.

Denote $M^{p,0}_{G}$ be the space with appropriate dimension of elementary process, $\theta_{t}=\sum_{i=0}^{n-1}\theta_{t_{i}}\mathbf{1}_{[t_{i},t_{i+1})}(t)$ with  each component of $\theta_{t_{i}}$ being in $L^{p}_{G}(\Omega_{t_{i}})$. Define norm
\begin{equation*}
\parallel \theta\parallel_{M^{p}_{G}}^{p}=\mathbb{E}[\int_{0}^{T}\vert \theta_{t}\vert^{p}dt], \theta\in M^{p,0}_{G},
\end{equation*}
and let $M^{p}_{G}$ denote the closure of $M^{p,0}_{G}$ under the norm $\parallel \cdot\parallel_{M^{p}_{G}}^{p}$.

In the following, denote $M^{p}_{G}(\mathbb{R}^{d})$(respectively, $M^{p}_{G}(\mathbb{R}^{d\times n})$ and  $M^{p}_{G}(\mathbb{R}^{n\times d\times d})$) the complete normed space under the norm
$\parallel \cdot\parallel_{M^{p}_{G}}$ with $\mathbb{R}^{d}$(respectively, $\mathbb{R}^{d\times n}$ and  $\mathbb{R}^{n\times d\times d}$)-valued processes.

For $\beta>0$ and $\eta\in M^{2}_{G}$, we let $M^{2,\beta}_{G}$ denote the space $M^{2}_{G}$ endowed
with the norm
\begin{equation*}
\parallel \eta\parallel_{M^{2,\beta}_{G}}:=\left\{\mathbb{E}[\int_{0}^{T}e^{\beta t}\vert\eta_{t}\vert^{2}dt]\right\}^{\frac{1}{2}}.
\end{equation*}

It is proved in Denis, Hu and Peng [2] 
that there exists
a weakly compact family $\mathcal{P}$ of probability measures defined on $(\Omega, \mathcal{B}(\Omega))$
such that
\begin{equation*}
\mathbb{E}[X]=\sup_{P\in\mathcal{P}}E_{P}[X], \text{ for } X\in L_{ip}(\Omega)=\bigcup_{n=1}^{\infty}L_{ip}(\Omega_{n}).
\end{equation*}
The natural choquet capacity is defined as
\begin{equation*}
c(A):=\sup_{P\in\mathcal{P}}P(A), \text{ for } A\in \mathcal{B}(\Omega).
\end{equation*}

\begin{defi} A set A is polar if
$c(A)=0$ and a property holds ``quasi-surely'' (q.s.) if it holds outside a polar set.
\end{defi}

Let's denote $\underline{\sigma}^{2}_{min}=\min_{1\leq i\leq d}\underline{\sigma}^{2}_{i}$, and $\bar{\sigma}^{2}_{max}=\max_{1\leq i\leq d}\bar{\sigma}^{2}_{i}$.

\section{Some Apriori Estimates of G-BSDEs}


Consider the following G-BSDE
\begin{equation}\label{G-BSDE integral form}
Y_{t}=\xi+\int_{t}^{T}f(s,Y_{s},Z_{s},\eta_{s})ds+\int_{t}^{T}g(s,Y_{s},Z_{s},\eta_{s}):d\langle B\rangle_{s}-\int_{t}^{T}Z^{*}_{s} dB_{s}+\int_{t}^{T}G(\eta_{s})ds-\frac{1}{2}\int_{t}^{T}\eta_{s}:d\langle B\rangle_{s},
\end{equation}
or, equivalently,
\begin{equation}\label{G-BSDE differential form}
-dY_{t}=f(t,Y_{t},Z_{t},\eta_{t})dt+g(t,Y_{t},Z_{t},\eta_{t}):d\langle B\rangle_{t}-Z^{*}_{t} dB_{t}+G(\eta_{t})dt-\frac{1}{2}\eta_{t}:d\langle B\rangle_{t}, Y_{T}=\xi,
\end{equation}
where

The terminal value is an $\mathcal{F}_{T}$-measurable random variable, $\xi: \Omega\mapsto \mathbb{R}^{n}$.

The generator $f$ maps $\Omega\times\mathbb{R}^{+}\times\mathbb{R}^{n}\times\mathbb{R}^{d\times n}\times\mathbb{D}^{n\times d\times d}$ onto $\mathbb{R}^{n}$ and is $\mathcal{B}\otimes\mathcal{B}^{n}\otimes\mathcal{B}^{d\times n}\otimes\mathcal{B}^{n\times d\times d}$-measurable.

The generator $g$ maps $\Omega\times\mathbb{R}^{+}\times\mathbb{R}^{n}\times\mathbb{R}^{d\times n}\times\mathbb{D}^{n\times d\times d}$ onto $\mathbb{D}^{n\times d\times d}$ and is $\mathcal{B}\otimes\mathcal{B}^{n}\otimes\mathcal{B}^{d\times n}\otimes\mathcal{B}^{n\times d\times d}$-measurable.

Here $\mathcal{B}$ is the $\sigma$-field of predictable sets of $\Omega\times[0,T]$.

Suppose that $\xi\in L^{2}_{G}(\Omega_{T})$, $f(\cdot,y,z,\eta),g(\cdot,y,z,\eta)\in M^{2}_{G}$ for each $y\in \mathbb{R}^{n}, z\in \mathbb{R}^{d\times n}, \eta\in\mathbb{D}^{n\times d\times d}$, and $f,g$ are uniformly Lipschitz; i.e., there exists $C>0$ such that for every $t$
\begin{equation*}
\begin{split}
&\vert f(\omega,t,y_{1},z_{1},\eta_{1})-f(\omega,t,y_{2},z_{2},\eta_{2})\vert\leq C(\vert y_{1}-y_{2}\vert+\vert z_{1}-z_{2}\vert+\vert \eta_{1}-\eta_{2}\vert), \forall (y_{1},z_{1},\eta_{1}), \forall (y_{2},z_{2},\eta_{2}),\\
&\vert g(\omega,t,y_{1},z_{1},\eta_{1})-g(\omega,t,y_{2},z_{2},\eta_{2})\vert\leq C(\vert y_{1}-y_{2}\vert+\vert z_{1}-z_{2}\vert+\vert \eta_{1}-\eta_{2}\vert), \forall (y_{1},z_{1},\eta_{1}), \forall (y_{2},z_{2},\eta_{2}).
\end{split}
\end{equation*}
Then we say $(\xi, f,g)$ are standard parameters for the G-BSDEs.

\begin{lemm}
$\forall \eta\in M^{1}_{G}(\mathbb{D}^{n\times d\times d})$, $\forall 0\leq t\leq s\leq T$, we have
\begin{equation}
\label{eta inequality 1}
\vert\int_{t}^{s}\eta_{r}:d\langle B\rangle_{r}\vert\leq K\int_{t}^{s}\vert\eta_{r}\vert dr,\\
\end{equation}
\begin{equation}
\label{eta inequality 2}
\int_{t}^{s}(\eta_{r})^{+}:\underline{\sigma}^{2}dr-\int_{t}^{s}(\eta_{r})^{-}:\bar{\sigma}^{2}dr\leq\int_{t}^{s}\eta_{r}:d\langle B\rangle_{r}\leq \int_{t}^{s}(\eta_{r})^{+}:\bar{\sigma}^{2}dr-\int_{t}^{s}(\eta_{r})^{-}:\underline{\sigma}^{2}dr,
\end{equation}
where $K$ is a constant, $\eta^{+}=((\eta^{1})^{+},\ldots,(\eta^{n})^{+})^{*}$ with $(\eta^{i})^{+}=diag((\eta^{i}_{1})^{+}, \ldots, (\eta^{i}_{d})^{+}), i=1,\ldots,n$, and
$\eta^{-}=((\eta^{1})^{-},\ldots,(\eta^{n})^{-})^{*}$ with $(\eta^{i})^{-}=diag((\eta^{i}_{1})^{-}, \ldots, (\eta^{i}_{d})^{-}), i=1,\ldots,n$.
\end{lemm}
Proof. Notice that the terms in (\ref{eta inequality 2}) are vectors, and the inequalities are for every elements of the vectors. It is easy to prove that (\ref{eta inequality 1}) and (\ref{eta inequality 2}) hold for $\eta\in M^{1,0}_{G}(\mathbb{D}^{n\times d\times d})$. Continuously extend them  to the case $\eta\in M^{1}_{G}(\mathbb{D}^{n\times d\times d})$ and we get  (\ref{eta inequality 1}) and (\ref{eta inequality 2}).

\begin{theo}\label{theorem beta(n) to infity}
$\theta, \zeta\in M^{2}_{G}$, $\theta_{s}, \zeta_{s}$ are continuous quasi-surely about $s$, and  $\parallel\theta\parallel_{M^{2}_{G}}, \parallel\zeta\parallel_{M^{2}_{G}}\neq 0$. Then there exists a sequence of $\beta(i)\rightarrow \infty$, as $i\rightarrow \infty$, such that
\begin{equation}\label{beta(n) to infity}
\lim_{i\rightarrow \infty}\frac{\mathbb{E}[\int_{t}^{T}e^{\beta(i)
s}\theta^{2}_{s}ds]}{\beta(i)\mathbb{E}[\int_{t}^{T}e^{\beta(i) s}\zeta^{2}_{s}ds]}=0.
\end{equation}
\end{theo}

Proof. $\forall \theta \in M^{2}_{G}$, there exists $\theta^{n}\xlongrightarrow[]{M^{2}_{G}} \theta, n\rightarrow \infty $,
with $\theta^{n}_{s}=\sum_{i=0}^{N_{n}-1}\theta^{n}_{s_{i}^{n}}\mathbf{1}_{[s_{i}^{n},s_{i+1}^{n})}(s)$, $\theta^{n}_{s_{i}^{n}}\in
\mathbb{L}^{2}(\Omega_{s_{i}^{n}})$, $s_{0}^{n}=t, s_{N_{n}}^{n}=T$. Then
\begin{equation*}
\lim_{n\rightarrow \infty}\mathbb{E}[\int_{t}^{T}e^{\beta s}(\theta^{n}_{s})^{2}ds]= \mathbb{E}[\int_{t}^{T}e^{\beta s}\theta_{s}^{2}ds],
\forall \beta>0.
\end{equation*}
And
\begin{equation*}
\begin{split}
&\mathbb{E}[\int_{t}^{T}e^{\beta s}(\theta^{n}_{s})^{2}ds]\\
=&\mathbb{E}[\int_{t}^{T}e^{\beta s}\sum_{i=0}^{N_{n}-1}(\theta^{n}_{s_{i}^{n}})^{2}\mathbf{1}_{[s_{i}^{n},s_{i+1}^{n})}(s)ds]\\
=&\mathbb{E}[\sum_{i=0}^{N_{n}-1}(\theta^{n}_{s_{i}^{n}})^{2}\int_{s_{i}^{n}}^{s_{i+1}^{n}}e^{\beta s}ds]\\
\leq& \max_{0\leq i\leq N^{n}-1}\mathbb{E}(\theta^{n}_{s_{i}^{n}})^{2}\int_{t}^{T}e^{\beta s}ds.
\end{split}
\end{equation*}
Let $C_{n}=\max\limits_{0\leq i\leq N^{n}-1}\mathbb{E}(\theta^{n}_{s_{i}^{n}})^{2}$. Since $\parallel\theta\parallel_{M^{2}_{G}}\neq 0$,  there exists $n$ large enough such that $C_{n}\neq 0$, and
\begin{equation*}
\mathbb{E}[\int_{t}^{T}e^{\beta s}(\theta^{n}_{s})^{2}ds]\leq C^{n}\int_{t}^{T}e^{\beta s}ds.
\end{equation*}

By the same reason, there exists $\bar{\zeta}^{n}\xlongrightarrow[]{M^{2}_{G}}  \zeta, n\rightarrow \infty$,
with $\bar{\zeta}^{n}_{s}=\sum_{i=0}^{M_{n}-1}\bar{\zeta}^{n}_{s_{i}^{n}}\mathbf{1}_{[s_{i}^{n},s_{i+1}^{n})}(s)$, $\bar{\zeta}^{n}_{s_{i}^{n}}\in
\mathbb{L}^{2}(\Omega_{s_{i}^{n}})$,  $s_{0}^{n}=t, s_{N_{n}}^{n}=T$, such that
\begin{equation*}
\lim_{n\rightarrow \infty}\mathbb{E}[\int_{t}^{T}e^{\beta s}(\bar{\zeta}^{n}_{s})^{2}ds]= \mathbb{E}[\int_{t}^{T}e^{\beta s}\zeta_{s}^{2}ds], \forall
\beta>0.
\end{equation*}
Let $(\zeta^{n}_{s})^{2}=\sum_{i=0}^{M_{n}-1}(\bar{\zeta}^{n}_{s_{i}^{n}})^{2}\mathbf{1}_{[s_{i}^{n},s_{i+1}^{n})}(s)+\frac{1}{n}$. Then
$(\zeta^{n}_{s})^{2}\geq \frac{1}{n}$, $n=1,2,\ldots$, and
\begin{equation*}
\lim_{n\rightarrow \infty}\mathbb{E}[\int_{t}^{T}e^{\beta s}(\zeta^{n}_{s})^{2}ds]= \mathbb{E}[\int_{t}^{T}e^{\beta s}\zeta_{s}^{2}ds], \forall
\beta>0.
\end{equation*}

Let $D_{n}=-\max\limits_{0\leq i\leq M^{n}-1}\{\mathbb{E}[-(\zeta^{n}_{s_{i}^{n}})^{2}]\}$.
Then $0<D_{n}<\infty$, and
\begin{equation*}
\begin{split}
&D_{n}\int_{t}^{T}e^{\beta s}ds-\mathbb{E}[\int_{t}^{T}e^{\beta s}(\zeta^{n}_{s})^{2}ds]\\
\leq&\mathbb{E}[D_{n}\sum_{i=0}^{M_{n}-1}\int_{s_{i}^{n}}^{s_{i+1}^{n}}e^{\beta s}ds-\sum_{i=0}^{M_{n}-1}(\zeta^{n}_{s_{i}^{n}})^{2}\int_{s_{i}^{n}}^{s_{i+1}^{n}}e^{\beta s}ds]\\
\leq&\sum_{i=0}^{M_{n}-1}[D_{n}+\mathbb{E}[-(\zeta^{n}_{s_{i}^{n}})^{2}]\int_{s_{i}^{n}}^{s_{i+1}^{n}}e^{\beta s}ds]\\
\leq& 0.
\end{split}
\end{equation*}
Hence
\begin{equation*}
\mathbb{E}[\int_{t}^{T}e^{\beta s}(\zeta^{n}_{s})^{2}ds]\geq D_{n}\int_{t}^{T}e^{\beta s}ds, n=1,2,\ldots.
\end{equation*}

Let
\begin{equation*}
B^{\beta}_{n}=\frac{\mathbb{E}[\int_{t}^{T}e^{\beta s}(\theta^{n}_{s})^{2}ds]}{\beta\mathbb{E}[\int_{t}^{T}e^{\beta
s}(\zeta^{n}_{s})^{2}ds]}.
\end{equation*}
Then
\begin{equation*}
B^{\beta}_{n}\leq \frac{C_{n}\int_{t}^{T}e^{\beta s}ds}{\beta D_{n}\int_{t}^{T}e^{\beta s}ds}=\frac{C_{n}}{\beta D_{n}}.
\end{equation*}
For any $n\in \mathbb{N}$ such that $C_{n}\neq 0$, let
\begin{equation}\label{beta(n) equal}
\beta(n)=n\frac{C_{n}}{D_{n}}>0,
\end{equation}
then
\begin{equation}\label{B_{n}}
B_{n}=\frac{\mathbb{E}[\int_{t}^{T}e^{\beta(n) s}(\theta^{n}_{s})^{2}ds]}{\beta(n)\mathbb{E}[\int_{t}^{T}e^{\beta(n)
s}(\zeta^{n}_{s})^{2}ds]}\leq \frac{1}{n} .
\end{equation}

Denote
\begin{equation*}
\begin{split}
T_{n}=&\frac{\mathbb{E}[\int_{t}^{T}e^{\beta(n) s}\theta^{2}_{s}ds]}{\beta(n)\mathbb{E}[\int_{t}^{T}e^{\beta(n)
s}\zeta^{2}_{s}ds]},\\
l_{n}=& \frac{\mathbb{E}[\int_{t}^{T}e^{\beta(n)
s}(\theta^{n}_{s})^{2}ds]}{\mathbb{E}[\int_{t}^{T}e^{\beta(n)
s}\theta^{2}_{s}ds]}, \\
m_{n}=& \frac{\mathbb{E}[\int_{t}^{T}e^{\beta(n)
s}(\zeta^{n}_{s})^{2}ds}{\mathbb{E}[\int_{t}^{T}e^{\beta(n) s}\zeta^{2}_{s}ds]}.
\end{split}
\end{equation*}
Then
\begin{equation}\label{T_{n}}
\begin{split}
T_{n}=\frac{m_{n}}{l_{n}}B_{n}.
\end{split}
\end{equation}

We say $l_{n}, n=1,2,\ldots$ is bounded. Otherwise there exists subsequence $l_{n_{i}}\rightarrow \infty, i\rightarrow \infty$, which means
$\frac{\mathbb{E}[\int_{t}^{T}e^{\beta(n_{i})
s}(\theta^{n_{i}}_{s})^{2}ds]}{\mathbb{E}[\int_{t}^{T}e^{\beta(n_{i})
s}\theta^{2}_{s}ds]}\rightarrow \infty, i\rightarrow \infty$.
Then for any $M>1$, there exist $I>0$, such that when $i>I$,
\begin{equation*}
\mathbb{E}[\int_{t}^{T}e^{\beta(n_{i})
s}(\theta^{n_{i}}_{s})^{2}ds]>M\mathbb{E}[\int_{t}^{T}e^{\beta(n_{i})
s}\theta^{2}_{s}ds]>0.
\end{equation*}

Let $\bar{\theta}^{n_{i}}_{s}=\sum_{i=0}^{N_{n}-1}\theta_{s_{i}^{n}}\mathbf{1}_{[s_{i}^{n},s_{i+1}^{n})}(s)$. Then $\bar{\theta}^{n_{i}}_{s}$
and $\theta^{n_{i}}_{s}$ are measurable on product measurable space $([t,T]\times\Omega, \mathcal{B}[t,T]\times\mathcal{F})$,  since $\mathbf{1}_{[s_{i}^{n},s_{i+1}^{n})}(s)$, $\theta_{s_{i}^{n}}$ and $\theta^{n}_{s_{i}^{n}}$ are measurable on $\mathcal{B}[t,T]\times\mathcal{F}$. Because $\theta_{s}$ is continuous quasi-surely about $s$, $\lim_{i\rightarrow \infty}\bar{\theta}^{n_{i}}_{s}(\omega)=\theta_{s}(\omega), q.s.$,  and  $\theta_{s}(\omega)$ is also $\mathcal{B}[t,T]\times\mathcal{F}$
measurable.

For $i>I$, let $E^{1}_{M,i}\times E^{2}_{M,i}=\{(s,\omega): \vert\theta^{n_{i}}_{s}(\omega)\vert>\sqrt{M}\vert\theta_{s}(\omega)\vert\}$. Then
\begin{equation*}
\mathbb{E}[\int_{t}^{T}(\theta^{n_{i}}_{s}-\theta_{s})^{2}ds]>(1-\frac{1}{\sqrt{M}})^{2}\mathbb{E}[\int_{E^{1}_{M,i}}\mathbf{1}_{E^{2}_{M,i}}(\theta^{n_{i}}_{s})^{2}ds]=(1-\frac{1}{\sqrt{M}})^{2}\mathbb{E}[\mathbf{1}_{E^{2}_{M,i}}\int_{E^{1}_{M,i}}(\theta^{n_{i}}_{s})^{2}ds].
\end{equation*}
While $\theta^{n}\xlongrightarrow[]{M^{2}_{G}} \theta, n\rightarrow \infty $, hence $\mathbb{E}[\mathbf{1}_{E^{2}_{M,i}}\int_{E^{1}_{M,i}}(\theta^{n_{i}}_{s})^{2}ds]\rightarrow 0$, $i\rightarrow \infty$, which means $\mu(E^{1}_{M,i})c(E^{2}_{M,i})\rightarrow 0$, $i\rightarrow \infty$, where $\mu$ is the Borel measure and $c$ is the Choquet capacity defined by $c(A)=\sup_{P\in\mathcal{P}}P(A)$, for  $A\in \mathcal{F}$. So
\begin{equation*}
\mathbb{E}[\int_{t}^{T}e^{\beta(n_{i})
s}(\theta^{n_{i}}_{s})^{2}ds]>M\mathbb{E}[\int_{t}^{T}e^{\beta(n_{i})
s}\theta^{2}_{s}ds], \text{ for all } i>I
\end{equation*}
is impossible.

Similarly, we can prove any convergent subsequence $l_{n_{i}}$, $\frac{1}{l_{n_{i}}}\nrightarrow\infty, i\rightarrow \infty$,
therefore $l_{n_{i}}\nrightarrow 0, i\rightarrow \infty$. And  $m_{n}$ is bounded by the same reason.

Since  $l_{n}, n=1,2,\ldots$ is bounded, there exists convergent subsequence. Let $n_{k_{i}}, i=1,2,\ldots$ be a subsequence of $n_{k}, k=1,2,\ldots$, such that $l_{n_{k_{i}}}\rightarrow a\neq 0, i\rightarrow \infty$. By (\ref{T_{n}}), we have $\lim_{i\rightarrow \infty}T_{n_{k_{i}}}=0$.

\begin{coro}\label{coro beta to infity} In the proof of theorem \ref{theorem beta(n) to infity}, $\beta(n)$ can be any real number such that
\begin{equation}\label{beta(n) larger}
\beta(n)\geq n\frac{C_{n}}{D_{n}}>0.
\end{equation}
Hence we have
\begin{equation}\label{beta to infity}
\lim_{\beta\rightarrow \infty}\frac{\mathbb{E}[\int_{t}^{T}e^{\beta
s}\theta^{2}_{s}ds]}{\beta\mathbb{E}[\int_{t}^{T}e^{\beta s}\zeta^{2}_{s}ds]}=0.
\end{equation}
\end{coro}

\begin{prop}\label{estimates of G-BSDE}
Let $((\xi^{i}, f^{i}, g^{i}); i=1,2)$ be two standard parameters of the G-BSDE (\ref{G-BSDE differential form}) and $(Y^{i}, Z^{i}, \eta^{i})$ be two solutions in space $M^{2}_{G}(\mathbb{R}^{n})\times M^{2}_{G}(\mathbb{R}^{d\times n})\times M^{2}_{G}(\mathbb{D}^{n\times d\times d})$ satisfying:

i) $Y^{i}_{t}, \eta^{i}_{t}, i=1,2$ are continuous in $t$ quasi-surely;

ii) If $Y^{1}=Y^{2}$, $t-a.e., \omega-q.s.$, then $\eta^{1}=\eta^{2}$, $t-a.e., \omega-q.s.$.

Put $\delta Y_{t}=Y^{1}_{t}-Y^{2}_{t}$, $\delta Z_{t}=Z^{1}_{t}-Z^{2}_{t}$, $\delta \eta_{t}=\eta^{1}_{t}-\eta^{2}_{t}$, $\delta f_{t}=f^{1}(t, Y^{1}_{t}, Z^{1}_{t},\eta^{1}_{t})-f^{2}(t, Y^{2}_{t}, Z^{2}_{t},\eta^{2}_{t})$, and $\delta g_{t}=g^{1}(t, Y^{1}_{t}, Z^{1}_{t},\eta^{1}_{t})-g^{2}(t, Y^{2}_{t}, Z^{2}_{t},\eta^{2}_{t})$. There exist $\beta_{0}(\delta Y, \delta \eta)$ such that when $\beta\geq\beta_{0}(\delta Y, \delta \eta)$, it follows that
\begin{equation}\label{estimate of Y}
\parallel \delta Y\parallel^{2}_{M^{2,\beta}_{G}}\leq \frac{1}{\underline{\sigma}^{2}_{min}}\left[ e^{\beta T}\mathbb{E}\vert \delta Y_{T}\vert^{2}+\frac{1}{\mu^{2}}\parallel \delta f\parallel^{2}_{M^{2,\beta}_{G}}+\frac{\bar{\sigma}^{2}_{max}}{\nu^{2}}\parallel \delta g\parallel^{2}_{M^{2,\beta}_{G}}\right],
\end{equation}
\begin{equation}\label{estimate of Z}
\parallel \delta Z\parallel^{2}_{M^{2,\beta}_{G}}\leq \frac{3}{\underline{\sigma}^{2}_{min}}\left[ e^{\beta T}\mathbb{E}\vert \delta Y_{T}\vert^{2}+\frac{1}{\mu^{2}}\parallel \delta f\parallel^{2}_{M^{2,\beta}_{G}}+\frac{\bar{\sigma}^{2}_{max}}{\nu^{2}}\parallel \delta g\parallel^{2}_{M^{2,\beta}_{G}}\right],
\end{equation}
\begin{equation}\label{estimate of eta}
\parallel \delta \eta\parallel^{2}_{M^{2,\beta}_{G}}\leq \frac{1}{\underline{\sigma}^{2}_{min}}\left[ e^{\beta T}\mathbb{E}\vert \delta Y_{T}\vert^{2}+\frac{1}{\mu^{2}}\parallel \delta f\parallel^{2}_{M^{2,\beta}_{G}}+\frac{\bar{\sigma}^{2}_{max}}{\nu^{2}}\parallel \delta g\parallel^{2}_{M^{2,\beta}_{G}}\right].
\end{equation}
\end{prop}

Proof. Let $(Y, Z, \eta)\in M^{2}_{G}(\mathbb{R}^{n})\times M^{2}_{G}(\mathbb{R}^{d\times n})\times M^{2}_{G}(\mathbb{D}^{n\times d\times d})$ be a solution of $(\ref{G-BSDE differential form})$. Then by (\ref{eta inequality 1}), there exists a constant $K>0$ such that
\begin{equation*}
\begin{split}
\vert Y_{t}\vert\leq& \vert \xi\vert+K\int_{0}^{T}\vert f(s,Y_{s},Z_{s},\eta_{s})\vert ds+K\int_{0}^{T} \vert g(s,Y_{s},Z_{s},\eta_{s})\vert ds\\
&+\sup_{0\leq t\leq T}\vert \int_{t}^{T}Z^{*}_{s} dB_{s}\vert+K\int_{0}^{T}\vert G(\eta_{s})\vert ds+K\int_{0}^{T} \vert\eta_{s}\vert ds
\end{split}
\end{equation*}
It follows from Burkholder-Davis-Gundy inequalities that there exists constants $0< k_{2}<K_{2}<\infty$ such that
\begin{equation*}
k_{2}\mathbb{E}\left[\int_{0}^{T}(Z_{s}Z^{*}_{s}):d\langle B\rangle_{s}\right]\leq \mathbb{E}\left[\sup_{0\leq t\leq T}\vert \int_{0}^{t}Z^{*}_{s} dB_{s}\vert^{2}\right]\leq K_{2}\mathbb{E}\left[\int_{0}^{T}(Z_{s}Z^{*}_{s}):d\langle B\rangle_{s}\right].
\end{equation*}
Since $\langle B\rangle_{t}$, $\underline{\sigma}^{2}$, and $\bar{\sigma}^{2}$ are diagonal matrices, only the diagonal elements works in the operation :,  so by (\ref{eta inequality 2}), we have
\begin{equation*}
k_{2}\mathbb{E}\left[\int_{0}^{T}(Z_{s}Z^{*}_{s}):\underline{\sigma}^{2}ds\right]\leq \mathbb{E}\left[\sup_{0\leq t\leq T}\vert \int_{0}^{t}Z^{*}_{s} dB_{s}\vert^{2}\right]\leq K_{2}\mathbb{E}\left[\int_{0}^{T}(Z_{s}Z^{*}_{s}):\bar{\sigma}^{2}ds\right].
\end{equation*}
Hence $\sup_{0\leq t\leq T}\vert \int_{0}^{t}Z^{*}_{s} dB_{s}\vert\in\mathbb{L}^{2}$, and $\int_{0}^{T}\vert G(\eta_{s})\vert ds, \int_{0}^{T}\vert \eta_{s}\vert ds\in L^{2}_{G}(\Omega_{T})$. Since $(\xi, f, g)$ are standard parameters, $\vert \xi\vert+\int_{0}^{T}\vert f(s,Y_{s},Z_{s},\eta_{s})\vert ds+\vert\int_{0}^{T} g(s,Y_{s},Z_{s},\eta_{s}):d\langle B\rangle_{s}\vert$ belongs to
$L^{2}_{G}(\Omega_{T})$ too. So we have $\sup_{0\leq t\leq T}\vert Y_{t}\vert\in \mathbb{L}^{2}$.

Applying It\^o's formula to $e^{\beta s}\vert \delta Y_{s}\vert^{2}$ (Li and Peng [7]), we have
\begin{equation}\label{Ito's formula to e delta Y}
\begin{split}
&e^{\beta t}\vert \delta Y_{t}\vert^{2}+\int_{t}^{T}\beta e^{\beta s}\vert \delta Y_{s}\vert^{2}ds +\int_{t}^{T} e^{\beta s}\delta Z_{s}\delta Z^{*}_{s}:d\langle B\rangle_{s}\\
=&e^{\beta T}\vert \delta Y_{T}\vert^{2}+\int_{t}^{T}2 e^{\beta s}\delta Y_{s}\cdot (f^{1}(s,Y^{1}_{s},Z^{1}_{s},\eta^{1}_{s})-f^{2}(s,Y^{2}_{s},Z^{2}_{s},\eta^{2}_{s}))ds\\
&+\int_{t}^{T}2 e^{\beta s}\delta Y_{s}\cdot (g^{1}(s,Y^{1}_{s},Z^{1}_{s},\eta^{1}_{s})-g^{2}(s,Y^{2}_{s},Z^{2}_{s},\eta^{2}_{s})):d\langle B\rangle_{s}-\int_{t}^{T}2 e^{\beta s}\delta Y^{*}_{s}\delta Z^{*}_{s}dB_{s}\\
&+\int_{t}^{T}2 e^{\beta s}\delta Y_{s}\cdot(G(\eta^{1}_{s})-G(\eta^{2}_{s}))ds-\int_{t}^{T}e^{\beta s}\delta Y_{s}\cdot\delta \eta_{s}:d\langle B\rangle_{s}.
\end{split}
\end{equation}
Since $\delta Y\in M^{2}_{G}(\mathbb{R}^{n}), \delta Z\in M^{2}_{G}(\mathbb{R}^{d\times n})$, $e^{\beta s}\delta Y^{*}_{s}\delta Z^{*}_{s}\in M^{1}_{G}$, the stochastic integral $\int_{t}^{T}2 e^{\beta s}\delta Y^{*}_{s}\delta Z^{*}_{s}dB_{s}$ is well defined.

If $\mathbb{E}[\int_{t}^{T}e^{\beta s}\vert \delta Y_{s}\vert^{2}ds]=0$, then $\int_{t}^{T}e^{\beta s}\vert \delta Y_{s}\vert^{2}ds=0$,
$\omega-q.s.$, and $\delta Y=0$, $t-a.e., \omega-q.s.$, by ii),  $\mathbb{E}[\int_{t}^{T}e^{\beta s}\vert \delta \eta_{s}\vert^{2}ds]=0$. And from (\ref{Ito's formula to e delta Y})
\begin{equation*}
\mathbb{E}[\int_{t}^{T} e^{\beta s}\delta Z_{s}\delta Z^{*}_{s}:d\langle B\rangle_{s}]\leq
\mathbb{E}[e^{\beta T}\vert \delta Y_{T}\vert^{2}].
\end{equation*}
Since
\begin{equation*}
\underline{\sigma}^{2}_{min}\int_{t}^{T} e^{\beta s} \vert\delta Z_{s}\vert^{2}ds\leq\int_{t}^{T} e^{\beta s}\delta Z_{s}\delta Z^{*}_{s}:d\langle B\rangle_{s},
\end{equation*}
we have (\ref{estimate of Z}).

If $\mathbb{E}[\int_{t}^{T}e^{\beta s}\vert \delta Y_{s}\vert^{2}ds]\neq0$, for given $\beta$, and $(Y^{i}, \eta^{i})_{i=1,2}$, there exists $C(\beta,(Y^{i}, \eta^{i})_{i=1,2})$ such that
\begin{equation}\label{C(beta)}
\begin{split}
&C(\beta,(Y^{i}, \eta^{i})_{i=1,2})\\
=&\frac{\mathbb{E}[\int_{t}^{T}2 e^{\beta s}\delta Y_{s}\cdot(G(\eta^{1}_{s})-G(\eta^{2}_{s}))ds-\int_{t}^{T}e^{\beta s}\delta Y_{s}\cdot\delta \eta_{s}:d\langle B\rangle_{s}]+\underline{\sigma}^{2}_{min}\mathbb{E}[\int_{t}^{T}e^{\beta s}\vert \delta \eta_{s}\vert^{2}ds]}{\mathbb{E}[\int_{t}^{T}e^{\beta s}\vert \delta Y_{s}\vert^{2}ds]},
\end{split}
\end{equation}
i.e.,
\begin{equation}\label{change}
\begin{split}
&\mathbb{E}[\int_{t}^{T}2 e^{\beta s}\delta Y_{s}\cdot(G(\eta^{1}_{s})-G(\eta^{2}_{s}))ds-\int_{t}^{T}e^{\beta s}\delta Y_{s}\cdot\delta \eta_{s}:d\langle B\rangle_{s}]\\
=&C(\beta,(Y^{i},\eta^{i})_{i=1,2})\mathbb{E}[\int_{t}^{T}e^{\beta s}\vert \delta Y_{s}\vert^{2}ds]-\underline{\sigma}^{2}_{min}\mathbb{E}[\int_{t}^{T}e^{\beta s}\vert \delta \eta_{s}\vert^{2}ds].
\end{split}
\end{equation}

From (\ref{Ito's formula to e delta Y}), we have for any constants $\mu, \nu$,
\begin{equation}
\begin{split}
&e^{\beta t}\vert \delta Y_{t}\vert^{2}+\int_{t}^{T}\beta e^{\beta s}\vert \delta Y_{s}\vert^{2}ds +\int_{t}^{T} e^{\beta s}\delta Z_{s}\delta Z^{*}_{s}:d\langle B\rangle_{s}\\
\leq&e^{\beta T}\vert \delta Y_{T}\vert^{2}+\mu^{2}\int_{t}^{T}e^{\beta s}\vert \delta Y_{s}\vert^{2}ds+\frac{1}{\mu^{2}}\int_{t}^{T}e^{\beta s}\vert \delta f_{s}\vert^{2}ds+\nu^{2}\bar{\sigma}^{2}_{max}\int_{t}^{T}e^{\beta s}\vert \delta Y_{s}\vert^{2}ds\\
&+\frac{\bar{\sigma}^{2}_{max}}{\nu^{2}}\int_{t}^{T}e^{\beta s}\vert \delta g_{s}\vert^{2}ds-\int_{t}^{T}2 e^{\beta s}\delta Y^{*}_{s}\delta Z^{*}_{s}dB_{s}\\
&+\int_{t}^{T}2 e^{\beta s}\delta Y_{s}\cdot(G(\eta^{1}_{s})-G(\eta^{2}_{s}))ds-\int_{t}^{T}e^{\beta s}\delta Y_{s}\cdot\delta \eta_{s}:d\langle B\rangle_{s},\\
\end{split}
\end{equation}
and further
\begin{equation}\label{deltaY ,delta eta, deltaZ}
\begin{split}
&e^{\beta t}\vert \delta Y_{t}\vert^{2}+(\beta-\mu^{2}-\nu^{2}\bar{\sigma}^{2}_{max})\int_{t}^{T} e^{\beta s}\vert \delta Y_{s}\vert^{2}ds +\int_{t}^{T} e^{\beta s}\delta Z_{s}\delta Z^{*}_{s}:d\langle B\rangle_{s}\\
\leq&e^{\beta T}\vert \delta Y_{T}\vert^{2}+\frac{1}{\mu^{2}}\int_{t}^{T}e^{\beta s}\vert \delta f_{s}\vert^{2}ds+\frac{\bar{\sigma}^{2}_{max}}{\nu^{2}}\int_{t}^{T}e^{\beta s}\vert \delta g_{s}\vert^{2}ds-\int_{t}^{T}2 e^{\beta s}\delta Y^{*}_{s}\delta Z^{*}_{s}dB_{s}\\
&+\int_{t}^{T}2 e^{\beta s}\delta Y_{s}\cdot(G(\eta^{1}_{s})-G(\eta^{2}_{s}))ds-\int_{t}^{T}e^{\beta s}\delta Y_{s}\cdot\delta \eta_{s}:d\langle B\rangle_{s}.\\
\end{split}
\end{equation}
Hence
\begin{equation}
\begin{split}
&(\beta-\mu^{2}-\nu^{2}\bar{\sigma}^{2}_{max})\mathbb{E}[\int_{t}^{T} e^{\beta s}\vert \delta Y_{s}\vert^{2}ds]\\
\leq&e^{\beta T}\mathbb{E}[\vert \delta Y_{T}\vert^{2}]+\frac{1}{\mu^{2}}\mathbb{E}[\int_{t}^{T}e^{\beta s}\vert \delta f_{s}\vert^{2}ds]+\frac{\bar{\sigma}^{2}_{max}}{\nu^{2}}\mathbb{E}[\int_{t}^{T}e^{\beta s}\vert \delta g_{s}\vert^{2}ds]\\
&+C(\beta,(Y^{i}, \eta^{i})_{i=1,2})\mathbb{E}[\int_{t}^{T}e^{\beta s}\vert \delta Y_{s}\vert^{2}ds]-\underline{\sigma}^{2}_{min}\mathbb{E}[\int_{t}^{T}e^{\beta s}\vert \delta \eta_{s}\vert^{2}ds].\\
\end{split}
\end{equation}
And then
\begin{equation}\label{deltaY ,delta eta}
\begin{split}
&(\beta-\mu^{2}-\nu^{2}\bar{\sigma}^{2}_{max}-C(\beta,(Y^{i}, \eta^{i})_{i=1,2}))\mathbb{E}[\int_{t}^{T} e^{\beta s}\vert \delta Y_{s}\vert^{2}ds]+\underline{\sigma}^{2}_{min}\mathbb{E}[\int_{t}^{T}e^{\beta s}\vert \delta \eta_{s}\vert^{2}ds]\\
\leq&e^{\beta T}\mathbb{E}[\vert \delta Y_{T}\vert^{2}]+\frac{1}{\mu^{2}}\mathbb{E}[\int_{t}^{T}e^{\beta s}\vert \delta f_{s}\vert^{2}ds]+\frac{\bar{\sigma}^{2}_{max}}{\nu^{2}}\mathbb{E}[\int_{t}^{T}e^{\beta s}\vert \delta g_{s}\vert^{2}ds].\\
\end{split}
\end{equation}

By (\ref{C(beta)}),
\begin{equation}
\begin{split}
&\vert C(\beta,(Y^{i}, \eta^{i})_{i=1,2})\vert\\
\leq&\frac{\frac{3}{2}\bar{\sigma}^{2}_{max}\mathbb{E}[\int_{t}^{T}e^{\beta s}\vert \delta Y_{s}\vert^{2}ds]+\frac{5}{2}\bar{\sigma}^{2}_{max}\mathbb{E}[\int_{t}^{T}e^{\beta s}\vert \delta \eta_{s}\vert^{2}ds]}{\mathbb{E}[\int_{t}^{T}e^{\beta s}\vert \delta Y_{s}\vert^{2}ds]}=\frac{3\bar{\sigma}^{2}_{max}}{2}+\frac{5\bar{\sigma}^{2}_{max}}{2}\frac{\mathbb{E}[\int_{t}^{T}e^{\beta s}\vert \delta \eta_{s}\vert^{2}ds]}{\mathbb{E}[\int_{t}^{T}e^{\beta s}\vert \delta Y_{s}\vert^{2}ds]}.
\end{split}
\end{equation}
By corollary \ref{coro beta to infity},
\begin{equation*}
\lim_{\beta\rightarrow \infty}\frac{\mathbb{E}[\int_{t}^{T}e^{\beta
s}\delta \eta^{2}_{s}ds]}{\beta\mathbb{E}[\int_{t}^{T}e^{\beta s}\delta Y^{2}_{s}ds]}=0,
\end{equation*}
so we can always choose $\beta_{0}(\delta Y, \delta \eta)$ large enough such that when $\beta\geq\beta_{0}(\delta Y, \delta \eta)$,  for any given $\mu, \nu$, $\beta-\mu^{2}-\nu^{2}\bar{\sigma}^{2}_{max}-C(\beta,(Y^{i}, \eta^{i})_{i=1,2})\geq \underline{\sigma}^{2}_{min}$.
By (\ref{deltaY ,delta eta}), we get (\ref{estimate of Y}) and  (\ref{estimate of eta}).

By (\ref{deltaY ,delta eta, deltaZ}), if  $C(\beta,(Y^{i}, \eta^{i})_{i=1,2})\leq0$, obviously we have (\ref{estimate of Z}).

Otherwise, it is easy to test that
\begin{equation*}
\begin{split}
&\int_{t}^{T}2 e^{\beta s}\delta Y_{s}\cdot(G(\eta^{1}_{s})-G(\eta^{2}_{s}))ds-\int_{t}^{T}e^{\beta s}\delta Y_{s}\cdot\delta \eta_{s}:d\langle B\rangle_{s}\\
\leq& 5\frac{\bar{\sigma}^{4}_{max}}{\underline{\sigma}^{2}_{min}}\int_{t}^{T} e^{\beta s}\vert\delta Y_{s}\vert^{2} ds+2\underline{\sigma}^{2}_{min}\int_{t}^{T} e^{\beta s}\vert\delta \eta_{s}\vert^{2} ds,
\end{split}
\end{equation*}
so by (\ref{deltaY ,delta eta, deltaZ}),  we have
\begin{equation}\label{deltaY ,delta eta, deltaZ for deltaZ 1}
\begin{split}
&e^{\beta t}\vert \delta Y_{t}\vert^{2}+(\beta-\mu^{2}-\nu^{2}\bar{\sigma}^{2}_{max})\int_{t}^{T} e^{\beta s}\vert \delta Y_{s}\vert^{2}ds +\int_{t}^{T} e^{\beta s}\delta Z_{s}\delta Z^{*}_{s}:d\langle B\rangle_{s}\\
\leq&e^{\beta T}\vert \delta Y_{T}\vert^{2}+\frac{1}{\mu^{2}}\int_{t}^{T}e^{\beta s}\vert \delta f_{s}\vert^{2}ds+\frac{\bar{\sigma}^{2}_{max}}{\nu^{2}}\int_{t}^{T}e^{\beta s}\vert \delta g_{s}\vert^{2}ds-\int_{t}^{T}2 e^{\beta s}\delta Y^{*}_{s}\delta Z^{*}_{s}dB_{s}\\
&+5\frac{\bar{\sigma}^{4}_{max}}{\underline{\sigma}^{2}_{min}}\int_{t}^{T} e^{\beta s}\vert\delta Y_{s}\vert^{2} ds+2\underline{\sigma}^{2}_{min}\int_{t}^{T} e^{\beta s}\vert\delta \eta_{s}\vert^{2} ds,
\end{split}
\end{equation}
and then
\begin{equation}\label{deltaY ,delta eta, deltaZ for deltaZ 2}
\begin{split}
&e^{\beta t}\vert \delta Y_{t}\vert^{2}+(\beta-\mu^{2}-\nu^{2}\bar{\sigma}^{2}_{max}-5\frac{\bar{\sigma}^{4}_{max}}{\underline{\sigma}^{2}_{min}})\int_{t}^{T} e^{\beta s}\vert \delta Y_{s}\vert^{2}ds +\int_{t}^{T} e^{\beta s}\delta Z_{s}\delta Z^{*}_{s}:d\langle B\rangle_{s}\\
\leq&e^{\beta T}\vert \delta Y_{T}\vert^{2}+\frac{1}{\mu^{2}}\int_{t}^{T}e^{\beta s}\vert \delta f_{s}\vert^{2}ds+\frac{\bar{\sigma}^{2}_{max}}{\nu^{2}}\int_{t}^{T}e^{\beta s}\vert \delta g_{s}\vert^{2}ds-\int_{t}^{T}2 e^{\beta s}\delta Y^{*}_{s}\delta Z^{*}_{s}dB_{s}\\
&+2\underline{\sigma}^{2}_{min}\int_{t}^{T} e^{\beta s}\vert\delta \eta_{s}\vert^{2} ds.
\end{split}
\end{equation}
We choose $\beta$ large enough such that $\beta-\mu^{2}-\nu^{2}\bar{\sigma}^{2}_{max}-5\frac{\bar{\sigma}^{4}_{max}}{\underline{\sigma}^{2}_{min}}>0$,
then we have (\ref{estimate of Z}).

\begin{rem} \end{rem}
By the proof of theorem \ref{estimates of G-BSDE}, we also have
\begin{equation*}
\begin{split}
\mathbb{E}[\sup_{0\leq t\leq T}\vert\delta Y_{t}\vert^{2}]\leq\mathbb{E}[\sup_{0\leq t\leq T}e^{\beta t}\vert\delta Y_{t}\vert^{2}]\leq 3\left[ e^{\beta T}\mathbb{E}(\vert \delta Y_{T}\vert)^{2}+\frac{1}{\mu^{2}}\parallel \delta f\parallel^{2}_{M^{2,\beta}_{G}}+\frac{\bar{\sigma}^{2}_{max}}{\nu^{2}}\parallel \delta g\parallel^{2}_{M^{2,\beta}_{G}}\right].
\end{split}
\end{equation*}

\begin{rem} \end{rem}
Let $\xi\in L_{ip}(\Omega_{T})$, and $M_{t}=\mathbb{E}_{t}[\xi]$.
by Peng [15], there exist $Z, \eta$ such that
\begin{equation}\label{G-martingale representation theorem lip}
M_{t}=\mathbb{E}[\xi]+\int_{0}^{t}(Z_{s})^{*} dB_{s}-[\int_{0}^{t}G(\eta_{s})ds-\frac{1}{2}\int_{0}^{t}\eta_{s}:d\langle B\rangle_{s}],
\end{equation}
where $M_{t}, \eta_{t}, i=1,2, \ldots$ are continuous in $t$ quasi-surely, and  if $M=0$, $t-a.e., \omega-q.s.$, then $\eta=0$, $t-a.e., \omega-q.s.$.

Then $(M,Z, \eta)$ is the solution of G-BSDE
\begin{equation}\label{G-BSDE with f=0 g=0}
M_{t}=\xi-\int_{t}^{T}(Z_{s})^{*} dB_{s}+\int_{t}^{T}G(\eta_{s})ds-\frac{1}{2}\int_{t}^{T}\eta_{s}:d\langle B\rangle_{s}
\end{equation}
with parameter $(\xi,0,0)$.

By proposition \ref{estimates of G-BSDE}, there exist $\beta_{0}(M, \eta)$ such that when $\beta\geq\beta_{0}(M, \eta)$
\begin{equation}\label{estimate of Y, Z, eta for Lip}
\parallel M\parallel^{2}_{M^{2,\beta}_{G}}+\parallel Z\parallel^{2}_{M^{2,\beta}_{G}}+\parallel \eta\parallel^{2}_{M^{2,\beta}_{G}}\leq \frac{5}{\underline{\sigma}^{2}_{min}}e^{\beta T}\mathbb{E}[\xi^{2}].
\end{equation}
Hence $\beta_{0}(M, \eta)$ and $M,Z, \eta$ are uniquely determined by  $\xi$,  and we also denote $\beta_{0}(M, \eta)$ as $\beta_{0}(\xi)$.

For any $\xi\in L^{2}_{G}(\Omega_{T})$, there exist $\xi^{n}\in L_{ip}(\Omega_{T})$ such that
$\xi^{n}\xlongrightarrow[]{L^{2}_{G}} \xi, n\rightarrow \infty$. For each pair $n,m$, there exist $\beta_{0}(\xi^{n},\xi^{m})$ such that when
$\beta\geq\beta_{0}(\xi^{n},\xi^{m})$, we have
\begin{equation}\label{estimate of Yn-Ym, Zn-Zm, etan-etam for Lip}
\parallel M^{n}-M^{m}\parallel^{2}_{M^{2,\beta}_{G}}+\parallel Z^{n}-Z^{m}\parallel^{2}_{M^{2,\beta}_{G}}+\parallel \eta^{n}-\eta^{m}\parallel^{2}_{M^{2,\beta}_{G}}\leq \frac{5}{\underline{\sigma}^{2}_{min}}e^{\beta T}\mathbb{E}[(\xi^{n}-\xi^{m})^{2}].
\end{equation}

Define
\begin{equation}\label{L 2, beta G}
\begin{split}
\mathcal{L}^{2}_{G}(\Omega_{T})=\{&\xi\in L^{2}_{G}(\Omega_{T}): \text{ there exist } \xi^{n} \in L_{ip}(\Omega_{T}) \text{ and } \beta<\infty \text{ such that }\\ &\xi^{n}\xlongrightarrow[]{L^{2}_{G}} \xi, n\rightarrow \infty \text{ and } \beta_{0}(\xi^{n},\xi^{m})\leq\beta, n,m=1,2,\ldots \}
\end{split}
\end{equation}

\section{G-martingale representation and existence and uniqueness of G-BSDEs under a strong condition}

\begin{theo}\label{G-martingale representation theorem for mathcal{L} 2}
For any $\xi\in \mathcal{L}^{2}_{G}(\Omega_{T})$, denote $M_{t}=\mathbb{E}_{t}[\xi]$, then there exist unique $(Z, \eta)\in  M^{2}_{G}(\mathbb{R}^{d})\times M^{2}_{G}(\mathbb{D}^{d\times d})$ such that
\begin{equation}\label{G-martingale representation for mathcal{L} 2}
M_{t}=M_{0}+\int_{0}^{t}Z^{*}_{s} dB_{s}-[\int_{t}^{T}G(\eta_{s})ds-\frac{1}{2}\int_{t}^{T}\eta_{s}:d\langle B\rangle_{s}],
\end{equation}
\end{theo}
Proof. For any $\xi\in \mathcal{L}^{2}_{G}(\Omega_{T})$, there exist $\xi^{n}\in L_{ip}(\Omega_{T})$ and $\beta<\infty$  such that $\xi^{n}\xlongrightarrow[]{L^{2}_{G}} \xi, n\rightarrow \infty$  and $\beta_{0}(\xi^{n},\xi^{m})\leq\beta, n=1,2,\ldots$. For every $\xi^{n}$,
by Peng [15], there exist $M^{n}, Z^{n}, \eta^{n}$ such that
\begin{equation}\label{G-martingale representation theorem lip n}
M^{n}_{t}=\mathbb{E}[\xi^{n}]+\int_{0}^{t}(Z^{n}_{s})^{*} dB_{s}-[\int_{0}^{t}G(\eta^{n}_{s})ds-\frac{1}{2}\int_{0}^{t}\eta^{n}_{s}:d\langle B\rangle_{s}],
\end{equation}
$M^{n}_{t}, \eta^{n}_{t}, i=1,2, \ldots$ are continuous in $t$ quasi-surely, and if $M^{m}=M^{n}$, $t-a.e., \omega-q.s.$, then $\eta^{m}=\eta^{n}$, $t-a.e., \omega-q.s.$.

Then $(M^{n},Z^{n}, \eta^{n})$ is the solution of G-BSDE
\begin{equation}\label{G-BSDE with f=0 g=0}
M^{n}_{t}=\xi^{n}-\int_{t}^{T}(Z^{n}_{s})^{*} dB_{s}+\int_{t}^{T}G(\eta^{n}_{s})ds-\frac{1}{2}\int_{t}^{T}\eta^{n}_{s}:d\langle B\rangle_{s}
\end{equation}
with parameter $(\xi^{n},0,0)$.

By proposition \ref{estimates of G-BSDE}, we have when $\bar{\beta}\geq \beta(\xi^{m},\xi^{n})$
\begin{equation*}
\parallel M^{m}-M^{n}\parallel^{2}_{M^{2,\bar{\beta}}_{G}}+\parallel Z^{m}-Z^{n}\parallel^{2}_{M^{2,\bar{\beta}}_{G}}+\parallel \eta^{m}-\eta^{n}\parallel^{2}_{M^{2,\bar{\beta}}_{G}}\leq e^{\bar{\beta} T}\frac{5}{\underline{\sigma}^{2}_{min}}\parallel \xi^{m}-\xi^{n}\parallel^{2}_{L^{2}_{G}},
\end{equation*}
and consequently,
\begin{equation}\label{constant beta}
\parallel M^{m}-M^{n}\parallel^{2}_{M^{2}_{G}}+\parallel Z^{m}-Z^{n}\parallel^{2}_{M^{2}_{G}}+\parallel \eta^{m}-\eta^{n}\parallel^{2}_{M^{2}_{G}}\leq e^{\bar{\beta} T}\frac{5}{\underline{\sigma}^{2}_{min}}\parallel \xi^{m}-\xi^{n}\parallel^{2}_{L^{2}_{G}}.
\end{equation}
Let $\bar{\beta}=\beta$, then (\ref{constant beta}) holds for the constant $\beta$ and $m,n=1,2,\ldots$, and $(M^{n},Z^{n}, \eta^{n})$ is a Cauchy sequence in $M^{2}_{G}$, so there exist $(M, Z, \eta)$ such that $(M^{n},Z^{n}, \eta^{n})\xlongrightarrow[]{M^{2}_{G}}(M, Z, \eta)$. Since
\begin{equation*}
\begin{split}
&\mathbb{E}\left(\int_{t}^{T}G(\eta^{n}_{s})ds-\int_{t}^{T}G(\eta_{s})ds\right)^{2}\\
\leq& \mathbb{E}\left(\int_{t}^{T}\vert G(\eta^{n}_{s})-G(\eta_{s})\vert ds\right)^{2}\\
\leq& \bar{\sigma}^{4}_{max} \mathbb{E}[\int_{t}^{T}\vert \eta^{n}_{s}-\eta_{s}\vert^{2}ds]\rightarrow 0,
\end{split}
\end{equation*}
we have $\int_{t}^{T}G(\eta^{n}_{s})ds \xlongrightarrow[]{\mathbb{L}^{2}_{G}}\int_{t}^{T}G(\eta_{s})ds$.
By Denis, Hu and Peng[2] proposition 17, there exists a subsequence  $\int_{t}^{T}G(\eta^{n_{i}}_{s})ds \rightarrow\int_{t}^{T}G(\eta_{s})ds, q.s.$.
Similarly, there exist a subsequence $\xi^{n_{i_{j}}}$ converging to  $\xi$ quasi-surely, a subsequence $\int_{t}^{T}(Z^{n_{i_{j_{k}}}}_{s})^{*} dB_{s}$ converging to  $\int_{t}^{T}Z^{*}_{s} dB_{s}$ quasi-surely, and a subsequence $\int_{t}^{T}\eta^{n_{i_{j_{k_{l}}}}}_{s}:d\langle B\rangle_{s}$ converging to  $\int_{t}^{T}\eta_{s}:d\langle B\rangle_{s}$ quasi-surely. For simplicity, we denote the index of the quasi-surely convergent subsequences as $k$. Then $M^{k}_{t}\rightarrow M_{t}, k\rightarrow\infty, q.s.$.

Since $(M^{n},Z^{n}, \eta^{n})$ satisfy (\ref{G-martingale representation for mathcal{L} 2}), $(M,Z, \eta)$ satisfy (\ref{G-martingale representation for mathcal{L} 2}), q.s..

The uniqueness of $(Z, \eta)$  follows from proposition \ref{estimates of G-BSDE}.

\begin{theo}\label{existence and uniqueness of G-BSDE}
Given standard parameters $(\xi,f,g)$, let $\Psi(y,z,\zeta)=\xi+\int_{t}^{T}f(s,y_{s},z_{s},\zeta_{s})ds+\int_{t}^{T}g(s,y_{s},z_{s},\zeta_{s}):d\langle B\rangle_{s}$.
Suppose for any $(y,z,\zeta)\in M^{2}_{G}(\mathbb{R}^{n})\times M^{2}_{G}(\mathbb{R}^{d\times n})\times M^{2}_{G}(\mathbb{D}^{n\times d\times d})$, $\Psi(y,z,\zeta)\in \mathcal{L}^{2}_{G}(\Omega_{T})$, and there exists $\beta>0$, for any $(y,z,\zeta),(y',z',\zeta')\in M^{2}_{G}(\mathbb{R}^{n})\times M^{2}_{G}(\mathbb{R}^{d\times n})\times M^{2}_{G}(\mathbb{D}^{n\times d\times d})$, $\beta_{0}(\Psi(y,z,\zeta),\Psi(y',z',\zeta'))\leq\beta$, then there exists a unique triplet $(Y,Z,\eta)\in M^{2}_{G}(\mathbb{R}^{n})\times M^{2}_{G}(\mathbb{R}^{d\times n})\times M^{2}_{G}(\mathbb{D}^{n\times d\times d})$ which solves G-BSDE(\ref{G-BSDE integral form}) in the sense of  $\mathcal{P}$-q.s., and $Y$ is a $\mathcal{P}$-q.s. continuous process.
\end{theo}

Proof. Firstly, we prove there exists a mapping from $M^{2}_{G}(\mathbb{R}^{n})\times M^{2}_{G}(\mathbb{R}^{d\times n})\times M^{2}_{G}(\mathbb{D}^{n\times d\times d})$ into $M^{2}_{G}(\mathbb{R}^{n})\times M^{2}_{G}(\mathbb{R}^{d\times n})\times M^{*}_{G}(\mathbb{D}^{n\times d\times d})$,
\begin{equation*}
\Phi: (y,z,\zeta)\rightarrow (Y,Z,\eta)
\end{equation*}
where $(Y,Z,\eta)$ is the solution of the G-BSDE(\ref{G-BSDE integral form}) with generator $f(t,y_{t},z_{t},\zeta_{t})$,  $g(t,y_{t},z_{t},\zeta_{t})$, i.e.,
\begin{equation*}
Y_{t}=\xi+\int_{t}^{T}f(s,y_{s},z_{s},\zeta_{s})dt+\int_{t}^{T}g(s,y_{s},z_{s},\zeta_{s}):d\langle B\rangle_{s}-\int_{t}^{T}Z^{*}_{s}dB_{s}+\int_{t}^{T}G(\eta_{s})ds-\frac{1}{2}\int_{t}^{T}\eta_{s}:d\langle B\rangle_{s}.
\end{equation*}

Since $\xi+\int_{0}^{T} f(s,y_{s},z_{s},\zeta_{s})dt+\int_{0}^{T}g(s,y_{s},z_{s},\zeta_{s}):d\langle B\rangle_{s}\in \mathcal{L}^{2}_{G}\subseteq L^{2}_{G}$,
we can define a G-martingale $M_{t}:=\mathbb{E}_{t}[\xi+\int_{0}^{T} f(s,y_{s},z_{s},\zeta_{s})dt+\int_{0}^{T}g(s,y_{s},z_{s},\zeta_{s}):d\langle B\rangle_{s}]$. By theorem \ref{G-martingale representation theorem for mathcal{L} 2}, there exist unique
$(Z,\eta)\in M^{2}_{G}(\mathbb{R}^{d\times n})\times M^{2}_{G}(\mathbb{D}^{n\times d\times d})$ such that
\begin{equation*}
M_{t}=M_{0}+\int_{0}^{t}Z^{*}_{s} dB_{s}-\int_{0}^{t}G(\eta_{s})ds+\frac{1}{2}\int_{0}^{t}\eta_{s}:d\langle B\rangle_{s}, \mathcal{P}-q.s.
\end{equation*}
Since for every $s\in[0,T]$, $\langle B\rangle_{s}$ is a diagonal matrix, only the diagonal elements enter the operation :, so the uniqueness of  $\eta$ means the diagonal elements is uniquely determined. Hence we choose $\eta$ to be a diagonal matrix process.

Define the process $Y$ by
\begin{equation*}
Y_{t}=M_{t}-\int_{0}^{t}f(s,y_{s},z_{s},\eta_{s})ds-\int_{0}^{t}g(s,y_{s},z_{s},\eta_{s}):d\langle B\rangle_{s},
\end{equation*}
which is $\mathcal{P}$-q.s. continuous by Li and Peng [7]. 
And $Y$ is also given by
\begin{equation}\label{estimate of Y }
Y_{t}=\mathbb{E}_{t}[\xi+\int_{t}^{T} f(s,y_{s},z_{s},\zeta_{s})dt+\int_{t}^{T}g(s,y_{s},z_{s},\zeta_{s}):d\langle B\rangle_{s}].
\end{equation}
So
\begin{equation*}
\begin{split}
&Y_{t}+\int_{t}^{T}Z^{*}_{s} dB_{s}-\int_{t}^{T}G(\eta_{s})ds+\frac{1}{2}\int_{t}^{T}\eta_{s}:d\langle B\rangle_{s}\\
=&M_{0}+\int_{0}^{T}Z^{*}_{s} dB_{s}-\int_{0}^{T}G(\eta_{s})ds+\frac{1}{2}\int_{0}^{T}\eta_{s}:d\langle B\rangle_{s}-\int_{0}^{t}f(s,y_{s},z_{s},\eta_{s})ds-\int_{0}^{t}g(s,y_{s},z_{s},\eta_{s}):d\langle B\rangle_{s}\\
=&M_{T}-\int_{0}^{t}f(s,y_{s},z_{s},\eta_{s})ds-\int_{0}^{t}g(s,y_{s},z_{s},\eta_{s}):d\langle B\rangle_{s}\\
=&\xi+\int_{t}^{T}f(s,y_{s},z_{s},\eta_{s})ds+\int_{t}^{T}g(s,y_{s},z_{s},\eta_{s}):d\langle B\rangle_{s}, \mathcal{P}-q.s.,
\end{split}
\end{equation*}
which is
\begin{equation*}
Y_{t}=\xi+\int_{t}^{T}f(s,y_{s},z_{s},\zeta_{s})ds+\int_{t}^{T}g(s,y_{s},z_{s},\zeta_{s}):d\langle B\rangle_{s}-\int_{t}^{T}Z^{*}_{s} dB_{s}+\int_{t}^{T}G(\eta_{s})ds-\frac{1}{2}\int_{t}^{T}\eta_{s}:d\langle B\rangle_{s}.
\end{equation*}
By (\ref{estimate of Y }), we have $\sup_{0\leq t\leq T}\vert Y_{t}\vert\in \mathbb{L}^{2}$.

Let $(y^{1},z^{1},\zeta^{1})$, $(y^{2},z^{2},\zeta^{2})$ be two elements of $M^{2}_{G}(\mathbb{R}^{n})\times M^{2}_{G}(\mathbb{R}^{d\times n})\times M^{2}_{G}(\mathbb{D}^{n\times d\times d})$, and let  $(Y^{1},Z^{1},\eta^{1})$ and $(Y^{2},Z^{2},\eta^{2})$ be the associated solutions. Since $f(y,z,\zeta)$, $g(y,z,\zeta)$ do not contain  $Y,Z,\eta$, applying proposition \ref{estimates of G-BSDE},
\begin{equation*}
\begin{split}
&\parallel \delta Y\parallel^{2}_{M^{2,\beta}_{G}}+\parallel \delta Z\parallel^{2}_{M^{2,\beta}_{G}}+\parallel \delta\eta\parallel^{2}_{M^{2,\beta}_{G}}\\
\leq&\frac{5}{\underline{\sigma}^{2}_{min}\mu^{2}}\mathbb{E}\int_{0}^{T}e^{\beta s}\vert f(s,y^{1}_{s},z^{1}_{s},\zeta^{1}_{s})-f(s,y^{2}_{s},z^{2}_{s},\zeta^{2}_{s})\vert^{2}ds]\\
+&\frac{5}{\underline{\sigma}^{2}_{min}\nu^{2}}\mathbb{E}\int_{0}^{T}e^{\beta s}\vert g(s,y^{1}_{s},z^{1}_{s},\zeta^{1}_{s})-g(s,y^{2}_{s},z^{2}_{s},\zeta^{2}_{s})\vert^{2}ds].
\end{split}
\end{equation*}
Since $f,g$ is uniformly Lipschitz in $y,z, \zeta$,
\begin{equation*}
\parallel\delta Y\parallel^{2}_{M^{2,\beta}_{G}}+\parallel\delta Z\parallel^{2}_{M^{2,\beta}_{G}}+\parallel \delta\eta\parallel^{2}_{M^{2,\beta}_{G}}\leq\frac{5K}{\underline{\sigma}^{2}_{min}}\left(\frac{1}{\mu^{2}}+\frac{1}{\nu^{2}}\right)
[\parallel \delta y\parallel^{2}_{M^{2,\beta}_{G}}+\parallel \delta z\parallel^{2}_{M^{2,\beta}_{G}}+\parallel \delta\zeta\parallel^{2}_{M^{2,\beta}_{G}}],
\end{equation*}
where $K$ is a constant. By the proof of proposition \ref{estimates of G-BSDE}, we can choose $\mu,\nu$ large enough such that
\begin{equation*}
\frac{5K}{\underline{\sigma}^{2}_{min}}\left(\frac{1}{\mu^{2}}+\frac{1}{\nu^{2}}\right)<1,
\end{equation*}
and then the mapping $\Phi$ is a contraction from $M^{2}_{G}(\mathbb{R}^{n})\times M^{2}_{G}(\mathbb{R}^{d\times n})\times M^{2}_{G}(\mathbb{D}^{n\times d\times d})$ onto itself and there exists a fixed point, which is the unique
solution of the G-BSDE.

\section*{Acknowledgments}

This research was supported by Beijing Natural Science Foundation Grant 1112009.

\end{document}